\documentclass[12pt,leqno]{amsart}
\usepackage{amsmath,amssymb,amsfonts,amsthm}
\usepackage{eucal,graphicx}

\newtheorem{lem}{Lemma}
\newtheorem{cor}[lem]{Corollary}
\newtheorem{prop}[lem]{Proposition}
\newtheorem{thm}[lem]{Theorem}

\theoremstyle{remark}
\newtheorem{exam}{Example}%[section]
\newcommand \myexam[1]{\smallskip\begin{exam}[\emph{#1}]}

\renewcommand{\phi}{\varphi}

\newcommand\inv{^{-1}}
\newcommand\setm{\setminus}
\newcommand\chiz{\chi^\bbZ}
\newcommand\bbR{\mathbb{R}}
\newcommand\bbZ{\mathbb{Z}}
\newcommand\cH{\mathcal{H}}

\allowdisplaybreaks

\begin{document}

\title{Nonattacking Queens in a Rectangular Strip}

\author{Seth Chaiken}
\address{Computer Science Department\\
The University at Albany (SUNY)\\
Albany, NY 12222, U.S.A.}
\email{\tt sdc@cs.albany.edu}

\author{Christopher R.H.\ Hanusa}
\address{Department of Mathematical Sciences\\
Binghamton University (SUNY)\\
Binghamton, NY 13902-6000, U.S.A.\\
Current address:
Department of Mathematics\\
Queens College\\
65-30 Kissena Blvd.\\
Flushing, NY 11367}
\email{\tt chanusa@qc.cuny.edu}

\author{Thomas Zaslavsky}
\address{Department of Mathematical Sciences\\
Binghamton University (SUNY)\\
Binghamton, NY 13902-6000, U.S.A.}
\email{\tt zaslav@math.binghamton.edu}

\begin{abstract}
 The function that counts the number of ways to place nonattacking identical chess or fairy chess pieces in a rectangular strip of fixed height and variable width, as a function of the width, is a piecewise polynomial which is eventually a polynomial and whose behavior can be described in some detail.  We deduce this by converting the problem to one of counting lattice points outside an affinographic hyperplane arrangement, which Forge and Zaslavsky solved by means of weighted integral gain graphs.  
 We extend their work by developing both generating functions and a detailed analysis of deletion and contraction for weighted integral gain graphs.  
 For chess pieces we find the asymptotic probability that a random configuration is nonattacking, and we obtain exact counts of nonattacking configurations of small numbers of queens, bishops, knights, and nightriders.
 \end{abstract}



\thanks{Preprint of version published in Ann. Comb. 14 (2010) 419-441, DOI
10.1007/s00026-11-0068-7, online Feb. 15, 2011 \copyright Springer Basel AG 2011, submitted March 15, 2007}

\maketitle
\pagestyle{headings}

\tableofcontents

%%%%%%%%%%%%%%%%%%%%%%%%
\section{Nonattacking pieces}

A classic enumeration problem is that of counting the configurations of eight nonattacking queens on a chessboard, and more generally of $n$ nonattacking queens on an $n \times n$ board.  While thinking about how to represent this problem geometrically we realized that a similar problem, a generalization in which the height of the board is fixed but the board's width can vary, is amenable to treatment by a recent geometrical counting method of Forge and Zaslavsky \cite{SOA}.  

Here is our variant question:  Place $m$ nonattacking chess queens in a rectangular strip of fixed height $m$ and variable width $n$.  In how many ways can this be done?  
We show that, as $n$ increases, the answer $\nu_m(n)$ becomes a polynomial in $n$.  Indeed, the number of ways is a sum of more and more terms as $n$ increases from $0$, each term being zero until $n$ reaches a threshold after which the term is a product of linear factors.  There is one term, $n^m$, with largest degree, and this term is present for all nonnegative $n$; thus one may say that the number of nonattacking placements is a piecewise polynomial of degree $m$.  Furthermore, the signs of the terms are such that there is no cancellation between terms.  All this applies to any chess piece with centrally symmetric moves, placed one to a row; that is, to all but the pawn.  
Indeed it applies to any fairy chess piece with suitable moves.  A fairy chess piece is an imaginary chess piece with an arbitrary movement rule; many such pieces have been tried and some found exciting; we mention the queen with a knight's move added to its usual repertory, or the nightrider, which extends a knight's move indefinitely.  

Forge and Zaslavsky counted lattice points inside a hypercube of variable side length but not lying in any of a finite set of affinographic hyperplanes (hyperplanes determined by constancy of the difference of two coordinates), by converting the problem into one about weighted integral gain graphs (which will be defined later).  Their principal result is that the counting function is a piecewise polynomial with predictable leading term, and is eventually a polynomial.  This yields our basic theorem because, in the space of configurations of chess pieces, attack is an affinographic property.  

From knowing the moves of our piece we also get more detailed information about the counting function.  For pieces like the queen, bishop, and nightrider, we find the exact value of $m$ at which the function becomes a polynomial (Proposition \ref{P:exact}).  When there is one piece in each row, the magnitude of the second coefficient of the counting polynomial has a simple combinatorial interpretation; it equals half the maximum total number of squares attacked nonhorizontally if the pieces are placed on a wide board and they do not attack any of the same squares (this is a special case of Proposition \ref{P:secondterm}).  This second coefficient is important because it determines the asymptotic probability that a random placement of one piece in each row will be nonattacking (Corollary \ref{C:probability}).  

Besides this, we extend the general theory of weighted integral gain graphs to treat generating functions and the details of the method of deletion and contraction; then we apply the latter to get exact formulas and generating functions for the number of nonattacking configurations with fixed small height for four chess and fairy chess pieces: the queen, bishop, knight, and nightrider.  
We do not find any general formula for the counting function; that seems impractically difficult, on the order of finding the chromatic number of a graph.

%%%%%%%%%%%%%%%%%%%%%%%%
\section{Fundamentals}

Here is the exact setup in complete generality.  We have a rectangular board with $m$ rows and $n$ columns, $m$ being a fixed and $n$ a variable positive integer.  We think of the board as the set $[n] \times [m] \subseteq \bbR^2$, where $[n]:=\{1,2,\ldots,n\}$.  We have a chess piece $P$---a fairy chess piece, which means its moves are arbitrary (subject to two rules to be stated); the set of moves is $A_P \subseteq \bbZ^2$ and signifies that if $P$ is located at a position $(j,k) \in \bbZ^2$, then it attacks all the squares in $(j,k) + A_P := \{ (j,k)+(h,i) : (h,i)\in A_P \}$ that lie in the board.  (We say ``moves to'' as a synonym of ``attacks''; the actual moves, if different from the attacked squares, as with a pawn, are irrelevant to our problem.) 
We need not assume $(0,0) \in A_P$; but we do have two requirements for the moves.  First, only finitely many moves can have any one nonzero value of the vertical displacement.  Second, the set of horizontal moves, $\{ (h,0) \in A_P \}$, should be either finite or the infinite set $\bbZ \times \{0\}$.  
We shall also assume that $A_P$ is centrally symmetric, but this is just a technical convenience to ensure that if one piece (a copy of $P$) attacks another, the attack is reciprocal, and we can ensure it, if necessary, through replacing $A_P$ by $A_P \cup (-A_P)$.
Besides this, for each row $j$ we specify a number $q_j$ which is the exact number of pieces to occupy that row.  (The case $q_j>1$ is possible only if $P$ has a finite set of horizontal moves.  Then the ``move'' $(0,0)$ is important.  It prevents two pieces from occupying the same square, and if it is not in the move set, the number of pieces on a square is limited only by the number of pieces in the row.)  

Thus, we can describe the entire configuration of $q := q_1+q_2+\cdots+q_m$ pieces by $q$ variables $x_j^r$, where $1 \leq j \leq m$ and $1 \leq r \leq q_j$ and
 $$
 x_j^r = \text{ the column occupied by the $r$th piece in row $j$}.
 $$
 We call the vector $x := (x_j^r)$ a \emph{labelled configuration}.  
Since the pieces are identical but the setup treats them as distinguishable, two labelled configurations describe the same configuration of pieces if they have the same occupied positions, that is, the same multiset of values of the variables in each row.  
Thus, the number of nonattacking configurations, $\nu_m(n)$, is the number of labelled configurations, $\lambda_m(n)$, divided by the number of labellings:
\begin{equation}\label{E:convert}
\nu_m(n) = \frac{\lambda_m(n)}{q_1! q_2! \cdots q_m!} \, .
\end{equation}

Treating a labelled configuration as a point in $[n]^q \subseteq \bbR^q$, the rules of nonattack are simple.  A configuration is nonattacking precisely when: 
 \begin{enumerate}
 \item[]\begin{enumerate}
 \item[(N1)]  if $1 \leq i < j \leq m$, then $(x_j^s-x_i^r,j-i) \notin A_P$ for all $r \in [q_i]$ and $s \in [q_j]$, and
 \item[(N2)]  when $r<s$ in $[q_j]$, then $(x_j^s-x_j^r,0) \notin A_P$.
 \end{enumerate}
 \end{enumerate}

Rule (N1) is a finite list of requirements because $j-i$ is bounded by $m-1$.  Write $j-i = k \neq 0$.  Suppose there are $t_k$ moves with height $k$ and they are $(\mu_{k1},k), \ldots, (\mu_{kt_k},k)$.  (Possibly $t_k=0$.)  Then (N1) says
 $$
 x_{i+k}^s \neq x_i^r + \mu_{kl} 
 $$
 for $l \in [t_k]$ and $i\in[m-k]$.

Rule (N2) is also finite.  If the set of horizontal moves is finite, say $\pm(\mu_{01},0),\ldots,\pm(\mu_{0t_0},0)$, (N2) says 
 $$
 x_j^s \neq x_j^r \pm \mu_{0l} 
 $$
 for $l \in [t_0]$.  If it is infinite, then it is $\bbZ \times \{0\}$ so (N2) simply requires that all $q_j \leq 1$.

\begin{thm} \label{T:polynomial}
The number of nonattacking configurations of pieces having $q_j$ pieces in row $j$ for each $j \in [m]$ is a piecewise polynomial function of $n$ and is a polynomial if $n$ is sufficiently large.  The piecewise polynomial has degree $q = \sum q_j$ and leading coefficient $(q_1!q_2!\cdots q_m!)\inv$\,.  A sufficient condition for $n$ to be sufficiently large is
$$
n \geq (q-1) \cdot \max \{ |\mu| : (\mu,k) \in A_P,\ 0 \leq k < m \},
$$
(but we take $1 \leq k < m$ if $P$ has unbounded horizontal moves or if all $q_j \leq 1$).
\end{thm}

Later on, in Corollary \ref{C:detailed}, we get a more detailed description of this counting function.

For a piece with bounded horizontal moves, the total number of nonattacking configurations of $q$ pieces is a piecewise polynomial given by summing the values given by Theorem \ref{T:polynomial} for all weak compositions $q_1+\cdots+q_m$ of $q$.  Thus:

\begin{cor}
The number of nonattacking configurations of $q$ pieces with a finite number of horizontal moves is, for large $n$, a polynomial which has degree $q$ and leading coefficient
$$
\sum_{\substack{q_1+\cdots+q_m=q \\ q_i\geq0 }} \frac{1}{q_1!q_2!\cdots q_m!}\ .
$$
\end{cor}

The simplest case is that in which every row has one piece, i.e., all $q_i=1$ and $q=m$.  Then the factor $q_1!q_2!\cdots q_m! = 1$, $\nu_m(n) = \lambda_m(n)$, and the piecewise polynomial is monic.  We solved this case for four different pieces, the queen, bishop, knight, and nightrider (a fairy chess piece with indefinitely extended knight's move) for small numbers of rows; see the last section for the results.  The calculations quickly become too complicated for hand solution.  
Hanusa has prepared a pair of computer programs that manipulate weighted integral gain graphs to produce exact answers and generating functions; they are available on the Web \cite{WIGGprogram}.  The first phase, in Java, carries out operations on weighted integral gain graphs.  The second phase, in Maple, uses John Stembridge's symmetric function package SF \cite{SF} to produce generating functions.
The calculations uncovered an unexpected fact:  in the polynomial that gives the number of configurations for large $n$, the second coefficient has a simple combinatorial meaning.  See Proposition \ref{P:secondterm}.

%%%%%%%%%%%%%%
\section{Gain graphs, affinographic hyperplanes, and the integral chromatic function}\label{ggah}

An \emph{affinographic hyperplane} in $\bbR^q$ is a hyperplane whose equation has the form $x_j = x_i + c$; it is called \emph{integral} when $c$ is an integer.    An \emph{arrangement} of hyperplanes is a finite set, considered in conjunction with the way it decomposes $\bbR^q$ into $d$-dimensional pieces (called \emph{regions} of the arrangement).   

An \emph{integral gain graph} $\Phi$ consists of an underlying graph $\|\Phi\|$, with vertex set $V = \{v_1,\ldots,v_q\}$ and edge set $E$, and a \emph{gain function} $\phi$, which assigns an integer to each oriented edge $e$, with the rule that $\phi(e\inv) = -\phi(e)$ where $e\inv$ denotes the same edge $e$ with the opposite orientation.  There are two kinds of edge: a \emph{link} has two different endpoints, while a \emph{loop} has both ends at the same vertex.
The gain of a path $P$, $\phi(P)$, is the sum of the gains of the edges, which must be oriented in a consistent direction along the path.  For brevity we write $\mu v_{j} v_{j'}$ for an edge with endpoints $v_j,v_{j'}$ and gain $\mu$ in the direction from $v_j$ to $v_{j'}$.

In \cite{SOA} integral gain graphs were used to treat integral affinographic hyperplane arrangements.  
To each edge $e$, with endpoints $v_i$ and $v_j$, there is a corresponding affinographic hyperplane in $\bbR^q$ whose equation is $x_j = x_i + \phi(e)$ if we calculate the gain with $e$ oriented from $v_i$ to $v_j$.  This hyperplane does not depend on the orientation of $e$ because, taking the opposite orientation $e\inv$ from $v_j$ to $v_i$, one gets the same equation.  Thus every integral gain graph has a corresponding affinographic hyperplane arrangement, and conversely every integral affinographic hyperplane arrangement has a corresponding integral gain graph.

We define $\chiz_\Phi(n)$, the \emph{integral chromatic function} of $\Phi$, as the number of mappings $x:  V \to [n]$ such that $x_{j'} \neq x_{j} + \mu$ for each edge $\mu v_{j} v_{j'}$.  (We may think of a function $x$ as an integer lattice point in the hypercube $[n]^m$ that is not in any of the affinographic hyperplanes that correspond to the edges of $\Phi$.)  
Let 
$$
n_0(\Phi) := \max\{ \phi(P) : P \text{ is a path in } \Phi \},
$$
the largest gain of any path.  According to \cite[Corollary 3.2 and the following paragraph]{SOA}, $\chiz_\Phi(n)$ is a polynomial on the domain $n \geq n_0(\Phi)$ but not on $n \geq n_0(\Phi)-1$.

To calculate $\chiz_\Phi(n)$ by deleting and contracting edges, as we do in Section \ref{dc}, we need an added feature: vertex weights.  
A \emph{weighted integral gain graph} $(\Phi,h)$ (which we shall herein call a ``gain graph'' for short) is an integral gain graph $\Phi$ which also has an integer weight $h_i$ assigned to each vertex $v_i$.  The integral gain graphs introduced earlier, and any integral gain graphs without explicit weights, implicitly have weight $0$ at every vertex.  
Given a weighted integral gain graph, $\chiz_{(\Phi,h)}(n)$ denotes the number of functions $x:  V \to \bbZ$ such that $x_{j} \neq x_{i} + \phi(f)$ for each edge $f$ from $v_i$ to $v_j$, and $h_i < x_i \leq n$ for every vertex.  We call $\chiz_{(\Phi,h)}$ the \emph{integral chromatic function} of the gain graph.  

The deletion-contraction identity to be stated in Equation \eqref{E:dc} later in the paper requires definitions of deletion and contraction of an edge $e$ in a weighted integral gain graph.  
The definition of \emph{deletion} is that $e$ is removed from the graph with no other change; we write $(\Phi,h)\setm e$ for the deleted gain graph.  Contraction is more complicated.  We write $(\Phi,h)/e$ for the contracted gain graph; but for computations a simpler notation is $(\Phi'',h'') := (\Phi,h)/e$, in which the underlying graph is $(V'',E'')$ and the gain function is $\phi''$.

Before defining contraction we must define switching.  Write $V = \{v_1,\ldots,v_q\}$.  
A function $\eta: V \to \bbZ$ is called a \emph{switching function}.  \emph{Switching} by $\eta$ changes the gain graph $(\Phi,h)$, with gain function $\phi$ and weight function $h$, to $(\Phi^\eta,h^\eta)$ with the same underlying graph but with gain function $\phi^\eta(f) := \phi(f) - \eta_i + \eta_j$, where $f$ is an edge oriented from $v_i$ to $v_j$, and with weights $h^\eta_k := h_k + \eta_k$.  

Now fix an edge $e$ with endpoints $v_i$ and $v_j$ and suppose it has nonnegative gain $\phi(e)$ when oriented from $v_i$ to $v_j$.  
One may think of contraction as a two-step process wherein one modifies $v_i$ and its incident edges and then merges $v_i$ into $v_j$.  
To do this we define a switching function $\eta$ by $\eta_k := 0$ for $k \neq i$ and $\eta_i := \phi(e)$.  
To \emph{contract} $e$, 
\begin{enumerate}
\item first switch $(\Phi,h)$ by $\eta$, which changes the gain of $e$ to $0$ and that of each other edge $f$ to $\phi''(f) := \phi(f)-\phi(e)$ if $f$ is a link incident with $v_i$ (we consider $e,f$ to be oriented away from $v_i$) and to $\phi''(f) := \phi(f)$ otherwise; the weight of every vertex stays the same except for $h_i$, which changes to $h^\eta_i = h_i+\phi(e)$; 
\item then coalesce the endpoints of $e$ to form a new vertex $v_{ij}'' \in V''$ and delete $e$, retaining the switched gains of all edges $f \neq e$ and the weights of all vertices $v_k \neq v_i, v_j$; and set the weight $h_{ij}''$ of the contracted vertex $v_{ij}''$ equal to $\max(h^\eta_i,h^\eta_j) = \max(h_i+\phi(e),h_j)$.  
\end{enumerate}
(All this is from \cite[Section 2]{SOA}, except that we have updated \cite{SOA}'s rooted integral gain graphs to \cite{HPT}'s more versatile weighted integral gain graphs.)  
Parallel edges with the same gain can be combined into one, because they represent the same constraint on coloring.

%%%%%%%%%%%%%%
\section{Formulas for non-attacking configurations}\label{nac}

The particular integral gain graph $\Phi$ we need for the chess problem has $q$ vertices, one for each piece, labelled $v_{jl}$ for $1 \leq j \leq m$ and $l \in [q_j]$ (representing piece $(j,l)$, i.e., piece $l$ in row $j$).  
There is an edge called $\mu v_{jl} v_{j'l'}$ for each $l \in [q_j],\ l' \in [q_{j'}]$, each distinct $j$ and $j'=j+k$ in $[m]$, and each $\mu$ such that $(\mu ,k) \in A_P$ (representing a possible attack $(\mu,k)$ from piece $(j,l)$ to piece $(j',l')$).  
There is also an edge called $\mu v_{jl}v_{jl'}$ for each $l \neq l'$ in $[q_j]$ and each $\mu$ such that $(\mu ,0) \in A_P$ (i.e., for each potential attack $(\mu,0)$ from piece $l$ to piece $l'$ in row $j$).  
The number $\mu$ is the gain of the edge; that is, $\phi(\mu v_{jl}v_{jl'}) = \mu$.  These edges are directed but they always come in opposite pairs, $\mu v_{jl} v_{j'l'}$ and $(-\mu) v_{j'l'} v_{jl}$; one should consider each pair as a single undirected edge.  

Define hyperplanes in $\bbR^q$, with coordinates $x_j^r$, by the equations
\begin{equation}
\tag{H1}
x_{i+k}^s - x_i^r = \mu_{kl}   \quad\text{ for } k \in [m-1], i\in[m-k], l\in[t_k], r\in[q_i], s\in[q_{i+k}],
\end{equation}
and, if $A_P$ is horizontally finite, 
\begin{equation}
\tag{H2} 
x_j^s - x_j^r = \pm \mu_{0l}   \quad\text{ for } j\in[m], 1 \leq r < s \leq q_j, l\in[t_0] .
\end{equation}
The set $\cH$ of all these hyperplanes is the arrangement of affinographic hyperplanes that corresponds to the integral gain graph $\Phi$ of the preceding paragraph.  

\begin{proof}[Proof of Theorem \ref{T:polynomial}]
A nonattacking labelled configuration is a point $x \in [n]^q \setm \big(\bigcup \cH \big)$, and conversely.  The main result of \cite{SOA}, Theorem 3.1, gives a formula for the number of such points, our $\lambda_m(n)$.  The formula shows that $\lambda_m(n)$ is a monic polynomial in $n$ of degree equal to the dimension of the hypercube $[n]^q$, if $n$ is large.  
Moreover, \cite[Corollary 3.2]{SOA} shows how $\lambda_m(n)$ behaves when $n$ is not  large; it is $\lambda_m(n) = \chiz_\Phi(n)$ where 
\begin{equation}\label{E:piecewise}
\chiz_\Phi(n) = \sum_{r=1}^{r_0} (-1)^{q-d_r} c_r  \prod_{i=1}^{d_r} \max(0, n-n_{ri}) ,
\end{equation}
in which $r_0$ is a positive integer, the $n_{ri}$ are nonnegative integers, $d_1=q$ and the first term is $n^q$, the other $d_r < q$, and $c_r$ is a positive integer with $c_1=1$; that explains our earlier remark that $\lambda_m(n)$ is a sum of more and more terms, all of degree less than $q$, as $n$ increases.
(Corollary \ref{C:detailed}, below, shows furthermore that $\min_i n_{ri} = 0$.)  
This justifies the description of $\lambda_m(n)$ as given by finitely many polynomials of the same degree and leading coefficient over the whole range $n \geq 0$.  

One gets $\nu_m(n)$ immediately from $\lambda_m(n)$ via Equation \eqref{E:convert}.

To prove the sufficient condition for $n$ to be large enough we  employ the integral gain graph $\Phi$.  
We can easily bound $n_0(\Phi)$.  The maximum gain of any edge is the maximum abscissa magnitude $|\mu|$ of any $(\mu, k) \in A_P$ that could represent an attack on the $m$-row board, which means that $|k| < m$.  By symmetry it is sufficient to consider $0 \leq k <m$.  If all $q_j\leq 1$, no horizontal attack is possible so we can restrict $k$ to be positive.  As a path has no more than $q-1$ edges, all path gains are $\leq (q-1) \cdot \max |\mu|$.
\end{proof}

The strip interpretation of $n_0(\Phi)$ is that it is the greatest breadth that can be attained, within the available height of $m$ rows, by a sequence of moves that hits each row $j$ at most $q_j$ times.

We can find $n_0(\Phi)$ exactly for some kinds of piece, including queens, bishops, and nightriders.  We suppose the nonhorizontal moves are all contained in the upper and lower of the four quadrants formed by the lines $y = \pm \alpha x$, where $\alpha\inv$ is a positive integer $b$.  That means any move $(x,y)$ satisfies $|x|/|y| \leq b$ if $y \neq 0$.  ($P$ may have bounded or unbounded horizontal moves.)  For instance, $\alpha = b\inv = 1$ for a queen or bishop.  
We further suppose that $(x, \pm bx) \in A_P$ for every nonzero integer $x$.  Finally, we assume there is exactly one piece in each row, so $q=m$.

\begin{prop}\label{P:exact} 
With  $P$ as described and with all $q_j=1$, $\nu_m(n)$ is a polynomial of degree $m$ on the domain $n \geq n_0(\Phi) = b \lfloor\frac{m^2-2}{2}\rfloor$ but not on $n \geq n_0(\Phi)-1$.
\end{prop}

\begin{proof} 
The gain graph has vertex set $\{ v_1, \ldots, v_m\}$ and edges $(\pm bh)e_{i,i+h}$ for $i, i+h \in [m]$ and $h\neq 0$.  
We want the maximum gain of a path; we do that by characterizing the paths of largest gain.  Let $W$ be such a path.  Note that $V$ is ordered by subscript and that, by maximality, every edge in $W$ is of the form $|bh|e_{i,i+h}$.  Thus we need only consider such edges, omitting any others; we may then simplify  the notation to $e_{i,i+h}$.

There are three kinds of internal vertices of $W$.  At a \emph{zig}, both neighbors are higher vertices.  At a \emph{zag}, both neighbors are lower.  At a \emph{smooth} vertex, one neighbor is higher and the other is lower.  There are also the two end vertices of $W$ and there might be isolated vertices, not in $W$. 

If $v_k$ is a smooth vertex, it is in edges $e_{ik}, e_{kj}$ of $W$ with $i<k<j$ or $i>k>j$.  Replacing $e_{ik} e_{kj}$ by  $e_{ij}$ leaves the gain unchanged.  Thus, we may assume $W$ has no smooth vertices.

Suppose that after eliminating smooth vertices there is a vertex $v_k \notin V(W)$.  We can extend $W$ from an endpoint $v_i$ by the edge $e_{ik}$ with gain $b|k-i|$.  That increases the gain, contradicting the maximality of $W$.  Consequently, no $v_k$ can exist, so $V(W)=V$, and indeed no path with maximum gain can have a smooth vertex.

All zigs are lower than all zags.  If not, say $k<k'$, $v_k$ is a zag, and $v_{k'}$ a zig.  Thus $W$ contains paths $e_{ik} e_{kj}$ with $i,j<k$ and $e_{i'k'}e_{k'j'}$ with $i',j'>k'$.  Replace these paths with $e_{ik'} e_{k'j}$ and $e_{i'k} e_{kj'}$.  This increases the gain.  

There are now two cases depending on the parity of $m$.  The endpoints $v_i, v_j$ of $W$ are both zigs, or both are zags, if $m$ is odd, and one is of each kind if $m$ is even.

Suppose $m$ is even.  Extend $W$ by the edge $e_{ij}$.  Now we have a cycle $C$ in which every vertex is a zig or a zag.  All zigs have lower indices, so they are $v_1, \ldots , v_{r-1}, v_r$ where $r = m/2$; and the zags are $v_m, v_{m-1}, \ldots , v_{m+r+1}$.  It is easy to calculate that 
$\phi(C) = b \big[ 4\binom{r}{2} + m \big] = bm^2/2$.  We recover $W$ by removing the edge of least gain, which is $e_{r,r+1}$ of gain $b$; thus, $n_0(\Phi) = b(m^2-2)/2$.

If $m$ is odd, we may suppose $v_i, v_j \in L$.  If there is $v_k \in L$ with $i < k$, then interchanging $v_i$ with $v_k$ in $W$ increases the gain.  It follows that the zigs are $v_1, \ldots, v_{(m-3)/2}$, $v_i = v_{(m-1)/2}$, and $v_j = v_{(m+1)/2}$.  Therefore, 
$$
n_0(\Phi) = \phi(W) = b \left\{ \left[ 2\binom{\frac{m+1}{2}}{2} - 1\right] + 2\binom{\frac{m-1}{2}}{2} + (m-1) \right\} ,
$$
which simplifies to $b(m^2-3)/2$.
\end{proof}

Suppose the inverted slope $\alpha\inv$ is not an integer.  Then $n_0(\Phi) < \alpha\inv \lfloor\frac{m^2-1}{2}\rfloor$ because the gain of $e_{ij}$ is the greatest width attainable with height $j-i$, that is, $\lfloor \alpha\inv(j-i) \rfloor$.  Thus, the gain of an edge in a path of maximum gain is sometimes $< \alpha\inv(j-i)$.  (This unavoidable fact requires some proof, which we omit.)  A consequence is that $\nu_m(n)$ may be a polynomial on a domain including $n \geq \alpha\inv \lfloor \frac{m^2-4}{2} \rfloor - 1$, contrary to Proposition \ref{P:exact}.  Nevertheless, we obtain a more general bound than that of Theorem \ref{T:polynomial}.  It is no longer necessary to assume a piece appears in every row but we do assume any row has at most one piece; thus $q$, the number of pieces, is also the number of occupied rows.

\begin{prop}\label{P:improved} 
If all nonhorizontal moves are contained between the lines $y=\pm \alpha x$ where $\alpha > 0$, and if all $q_j \leq 1$, then $\nu_m(n)$ is a polynomial of degree $q$ on the range $n \geq \alpha\inv \lfloor\frac{m^2-2}{2}\rfloor$.
\end{prop}

The proof is similar enough to that of Proposition \ref{P:exact} that it can be omitted.

%%%%%%%%%%%%%%
\section{The method of deletion and contraction}\label{dc}

To calculate the integral chromatic function we need vertex weights and the formula for deletion and contraction of an edge in a weighted integral gain graph:
\begin{equation}\label{E:dc}
\chiz_{(\Phi,h)}(n) = \chiz_{(\Phi,h)\setm e}(n) - \chiz_{(\Phi,h)/e}(n) .
\end{equation}
This formula is valid for any link $e$, by \cite[Theorem 3.4]{SOA}, and translates into a formula for $\nu_m(n)$ by Equation \eqref{E:convert}.  

Another useful formula is multiplicativity.  If $\Phi$ is the disjoint union of $\Phi_1, \ldots, \Phi_r$, then
\[
\chiz_{(\Phi,h)} = \chiz_{(\Phi_1,h)} \cdots \chiz_{(\Phi_r,h)} .
\]

Now we begin to compute $\chiz_{(\Phi,h)}$.
In our calculations the notation $x^+$ means the positive part of the real number $x$, that is, 
$$
x^+ := \max(x,0).
$$  

The effect of a loop is simple.  If it has nonzero gain $k$, it can be discarded because it stands for a constraint $x_i \neq x_i + k$, which is always satisfied.  A loop with zero gain, however, forces $\chiz_\Phi(n) = 0$ for all $n$ because the corresponding constraint is unsatisfiable.

\myexam{A Vertex with Loops}\label{X:vertex}
If $\Phi$ has one vertex $v_1$, all edges are loops.  Thus, 
\[
\chiz_{(\Phi,h)}(n) = \begin{cases} (n-h_1)^+ &\text{ if there is no zero loop, } \\ 0 &\text{ if there is a zero loop.}\end{cases}
\]
\end{exam}

\myexam{Isolated Vertices with Loops}\label{X:isolated}
If $\Phi$ has no edges other than loops, then every component has a single vertex.  By Example \ref{X:vertex} and multiplicativity of the integral chromatic function,
\[
\chiz_{(\Phi,h)}(n) = \begin{cases} (n-h_1)^+ \cdots (n-h_q)^+ &\text{ if there are no zero loops, } \\ 0 &\text{ if there is a zero loop.}\end{cases}
\]
This gain graph is the basic type to which we reduce all other gain graphs by deletion and contraction.
\end{exam}

\myexam{A Multiple Edge}\label{X:multiple}
Suppose $\Phi$ has two vertices and all its edges have the form $e_i = \mu_iv_1v_2$, no two gains being equal.  Let $M = \{\mu_i\}$ be the set of gains.  If we contract one edge, the other edges become loops with nonzero gain, so may be ignored.  The switched weights are $(h_1+\mu_i,h_2)$ if $\mu_i \geq 0$ and $(h_1,h_2+|\mu_i|)$ if $\mu_i \leq 0$; the contracted weight is the greater of the two switched weights.  Thus, the deletion-contraction formula gives
\begin{align*}%\label{E:}
\chiz_{(\Phi,h)}(n) =\ &(n-h_1)^+ (n-h_2)^+ \\
&- \left\{ \sum_{i: \mu_i\geq0} \big(n-\max(h_1+\mu_i,h_2)\big)^+ 
+ \sum_{i: \mu_i<0} \big(n-\max(h_1,h_2+|\mu_i|)\big)^+ \right\} .
\end{align*}
If the original weights are both $0$, this simplifies to
\[
\chiz_{(\Phi,0)}(n) = (n^+)^2 - \sum_{i} (n-|\mu_i|)^+ .
\]
\end{exam}

\begin{figure}[hbt]
\includegraphics[scale=.8]{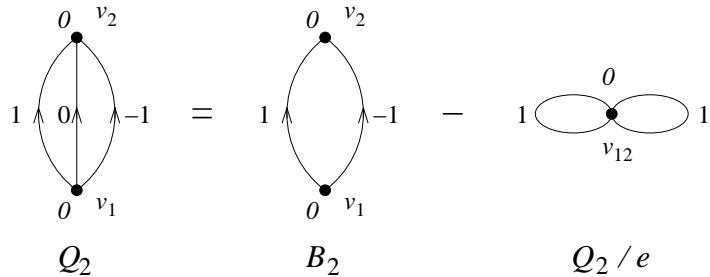}
\caption{The weighted integral gain graph $Q_2$, of order $2$ with edges having gains $-1,0,1$ from $v_1$ to $v_2$ and vertex weights all $0$, reduced by deletion and contraction of $e = 0v_1v_2$.  $B_2$ denotes $Q_2 \setm e$.}
\label{F:q2example}
\end{figure}

\begin{figure}[hbt]
\includegraphics[scale=.8]{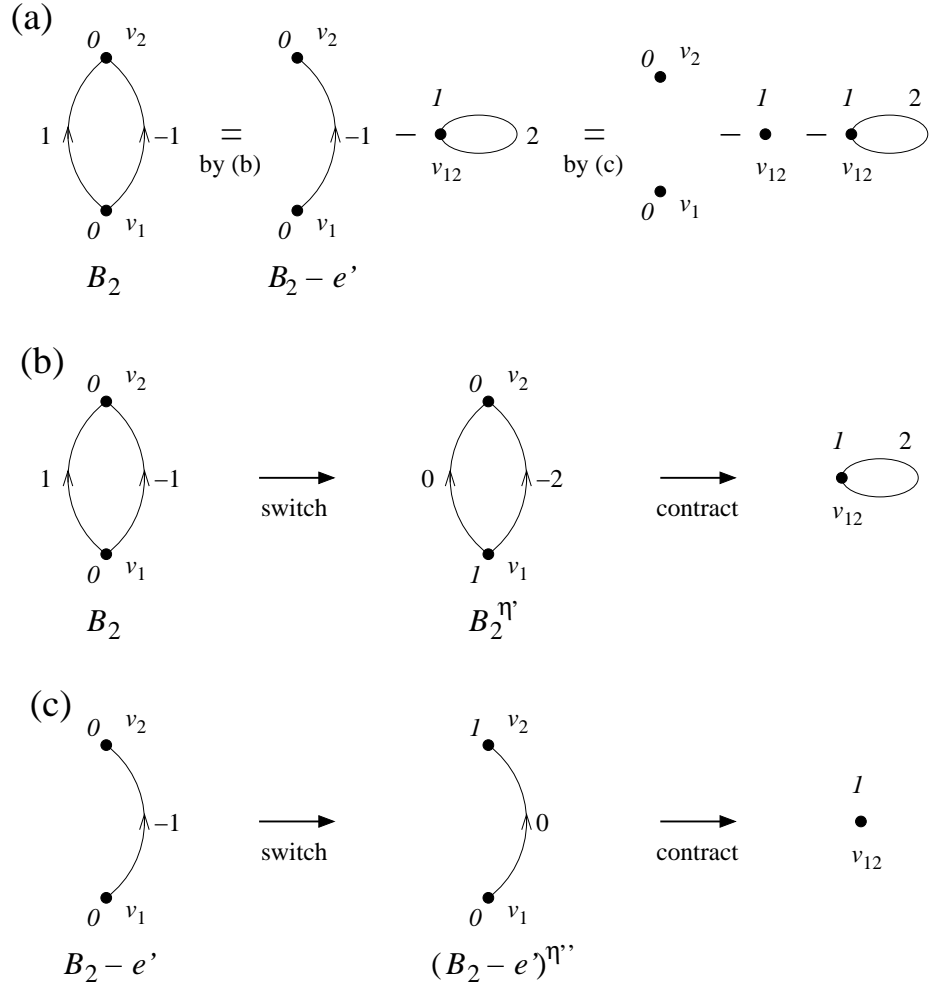}
\caption{(a)  The weighted integral gain graph $B_2$, of order $2$ with edges having gains $-1,1$ from $v_1$ to $v_2$ and zero weights, reduced by successive deletion and contraction of $e' = 1v_1v_2$, then $e'' = -1v_1v_2$.  
(b)  $B_2$ switched by $(\eta'_1,\eta'_2)=(1,0)$ so $e'$ can be contracted.  
(c)  The contraction process for $e''$, switching $B_2 \setm e'$ by $(\eta''_1,\eta''_2)=(0,1)$ so $e''$ has gain $0$, then contracting $e''$.}
\label{F:b2example}
\end{figure}

Figures \ref{F:q2example} and \ref{F:b2example} give examples of the deletion-contraction process symbolically, in the form of graph equations, and show the details of the contraction process with its two stages of switching (except for a zero edge) and contracting.  Symbolically, $Q_2$ is the difference between the deletion and the contraction of $e$, and so forth.

\myexam{A Small Problem of Queens and Bishops}\label{X:small}
To further explain the contraction-deletion method, we treat in detail the cases of queens and bishops with $m=2$ and $q_1=q_2=1$, that is, one piece in each of two rows, so that the (weighted integral) gain graph has order $2$.  
Since there is one piece per row, $\nu_m(n) = \chiz_\Phi(n)$.

The gain graph for each has vertices $v_1$ and $v_2$, one for each row, and all vertex weights equal to $0$.  The gain graph $Q_2$ for queens (see Figure \ref{F:q2example}) has edges $0v_1v_2$, $1v_1v_2$, $-1v_1v_2$ corresponding to the constraints $x_2 \neq x_1$, $x_2 \neq x_1+1$, $x_2 \neq x_1-1$ because a queen going up one row moves $0$, $1$, or $-1$ squares to the right.  
The gain graph $B_2$ for bishops (Figure \ref{F:b2example}(a)) has the same edges but without the one with zero gain since a bishop cannot move vertically.  From Example \ref{X:multiple} we get the formulas 
\begin{align} \label{E:small}
\nu_{B_2}(n) = \chiz_{B_2}(n) = n^2 - 2(n-1)^+ , \\
\nu_{Q_2}(n) = \chiz_{Q_2}(n) = n^2 - n - 2(n-1)^+ .
\end{align}

In the contraction-deletion process, the first step is to reduce the queens graph to the bishops graph by deleting and contracting $0v_1v_2$.  Figure \ref{F:q2example} shows this reduction.  The loops, which do not influence the integral chromatic function, can be deleted.  We deduce the formula
\[
\chiz_{Q_2}(n) = \chiz_{B_2}(n) - n ,
\]
because the contracted weight is $\max(0,0)$.

The next edge to go is $1v_1v_2$.  In order to contract it we have to switch by $\eta'_1 = 1$, $\eta'_2 = 0$, giving the switched graph $B_2^{\eta'}$ shown in Figure \ref{F:b2example}(b).  Then we can contract, as shown in the figure.  The weight in the contraction is $h_{12} = \max(1,0) = 1$ because the corresponding variable $x_{12}$ has to satisfy both bounds $x_1^{\eta'} > 1$ and $x_2^{\eta' }> 0$.  We deduce that
\[
\chiz_{B_2}(n) = \chiz_{-1v_1v_2}(n) - (n-1)^+ .
\]
The notation $-1v_1v_2$ here signifies a graph with two vertices and one edge whose gain is $-1$ from $v_1$ to $v_2$, or $+1$ in the opposite direction.  

We delete and contract $-1v_1v_2$ in Figure \ref{F:b2example}(a) and (c).   Switching by $\eta''_1 = 0, \eta''_2 = 1$ (since the positive direction of the edge is from $v_2$ to $v_1$) gives the gains and weights in Figure \ref{F:b2example}(c), whence comes the symbolic reduction formula in the figure (again, the loops should be ignored) and thence the formula 
\[
\chiz_{-2v_1v_2}(n) = n^2 - (n-1)^+ .
\]
Combining the steps we obtain \eqref{E:small}, as we ought.
\end{exam}

The contraction-deletion formula implies a slightly more precise description of the integral chromatic function, and thence of the counting function $\nu_m(n)$, than that stated in \cite[Corollary 3.2]{SOA}, quoted earlier at Equation \eqref{E:piecewise}.  
The improvement comes from giving a new proof of \eqref{E:piecewise} based on deletion and contraction.

\begin{cor}\label{C:detailed}
The integral chromatic function $\chiz_{(\Phi,0)}(n)$ of a nonempty $0$-weighted integral gain graph has the form of Equation \eqref{E:piecewise} where the $n_{ri}$ are nonnegative integers with all $n_{1i} = 0$.
\end{cor}

\begin{proof}
In \eqref{E:dc}, deletion preserves the order of the graph and the coefficient sign, while contraction changes the sign and the parity of the order.  Thus, after eliminating all links by repeated deletion and contraction, we have an alternating sum of integral chromatic functions $\chiz_{(\Psi_j,h_j)}$ of linkless integral gain graphs.  The sign of a term depends on the number of vertices lost in contraction.  Only one of the graphs has $q$ vertices; that is the one obtained by deleting every edge, and in it the weights remain $0$ as in the original gain graph.   It follows that
\begin{equation}\label{E:intermediate}
\chiz_{(\Phi,0)}(n) = \sum_j  (-1)^{q-d_j}  \chiz_{(\Psi_j,h_j)}(n) ,
\end{equation}
summed over all the linkless gain graphs that result from the deletion-contraction process, where $d_j$ is the order of $\Psi_j$.  Every one of these graphs has at least one vertex so $d_j>0$.  

The way weights contract ensures that the weights in each $\Psi_j$ are nonnegative, and also, since the original weights are all zero, that the minimum weight in $\Psi_j$ is $0$.  The integral chromatic function of a linkless graph $(\Psi,h)$ equals $\prod_{v \in V} (n-h_v)^+$ unless it is identically zero (see Example \ref{X:isolated}).  The corollary follows by substituting this expression in \eqref{E:intermediate} and collecting equal expressions.
\end{proof}

For sufficiently large integers $n$, the superscripts $+$ can be dropped.  When all the terms are multiplied out, there is no cancellation; the resulting formula is 
\begin{equation}\label{E:expanded}
\chiz_{(\Phi,0)}(n) =  \sum_{d=0}^q  (-1)^{q-d}  n^d  \sum_{j=1}^k  c_j \sigma_{d_j-d}(n_{j1},\ldots,n_{jd_j})  
 \text{ for } n \geq \max n_{ji} ,
\end{equation}
where $\sigma_{d_j-d}$ denotes the elementary symmetric function of degree ${d_j-d}$ (identically $0$ if $d<0$ or $d>d_j$, $1$ if $d=d_j$), all $c_j>0$, and all $n_{ji}\geq0$.

%%%%%%%%%%%%%%%%%%%%
\section{Generating functions}

As we have seen, the counting functions $\nu_m(n)$ that we get, and more generally those that arise in counting lattice points outside an integral affinographic hyperplane arrangement \cite{SOA}, are integer-weighted sums of terms of the form
\begin{equation}\label{E:term}
a_n = (n-n_1)^+ (n-n_2)^+ \cdots (n-n_r)^+ 
\end{equation}
where $n_1 \leq n_2 \leq \cdots \leq n_r$.  
The generating function of such a term has a simple form, 
${p(t)}/{(1-t)^{r+1}},$ 
where $p(t)$ is a polynomial which we can calculate exactly.   
That enables us to write down the generating function of $\nu_m$, which is  
\[
N_m(t) := \sum_{n=0}^\infty \nu_m(n) t^n = \frac{P(t)}{(1-t)^{q+1}} \, ,
\]
where $q$ is the number of vertices of the gain graph, i.e., the number of pieces on the board.  
This is a special case of the general behavior of the generating function of lattice-point counts of a weighted integral gain graph, $X_{(\Phi,h)}(t) := \sum_{n=0}^\infty \chiz_{(\Phi,h)}(n) t^n$, which has the form $X_{(\Phi,h)}(t) = P_{(\Psi,h)}(t) / (1-t)^{q+1}$ where the numerator is a polynomial.  The deletion-contraction formula \eqref{E:dc} extends immediately to generating functions: 
\begin{equation*}%\label{E:gfdc}
X_{(\Phi,h)}(t) = X_{(\Phi,h)\setm e}(t) - X_{(\Phi,h)/e}(t) 
\end{equation*}
for a link $e$; hence, the numerator polynomials satisfy
\[
P_{(\Psi,h)}(t) = P_{(\Psi,h)\setm e}(t) - (1-t) P_{(\Psi ,h) / e}(t) .
\]

Let us develop the generating function of a term of the form \eqref{E:term}.  
Let $d_i = n_r - n_i$ for $1 \leq i < r$ and let $s_k$ be the $k$th elementary symmetric function of $d_1, \ldots, d_{r-1}$.  Then 
\begin{equation}\label{E:plusgf}
\sum_{n=0}^\infty a_n t^n = t^{n_r} \sum_{j=1}^{r} (-1)^{r-j} s_{r-j} \, \frac{A_j(t)}{(1-t)^{j+1}}
\end{equation}
where $A_j(t)$ is the Eulerian polynomial, i.e., the polynomial such that 
$$\sum_{n=0}^\infty n^j t^n = \frac{A_j(t)}{(1-t)^{j+1}}$$ 
(see \cite{Comtet,EC1}).  The first few Eulerian polynomials are $A_0(t) = 1$, $A_1(t) = t$, $A_2(t) = t^2 + t$, $A_3(t) = t^3 + 4t^2 + t$, and $A_4(t) = t^4 + 11t^3 + 11t^2 + t$.

The proof of \eqref{E:plusgf} is by substituting for $n$ the variable $p := n - n_r$ and multiplying out the resulting product $p(p+d_1)\cdots(p+d_{r-1})$, in which the coefficient of $p^j$ is $s_{r-j}$.

The most useful cases for our examples are the generating function of $(n-n_1)^+$, which is
$$
\frac{t^{n_1+1}}{(1-t)^2},
$$
that of $n (n-n_2)^+$, which is
$$
\frac{t^{n_2+1} [ (n_2+1)t - (n_2-1) ]}{(1-t)^3},
$$
and that of $n^2 (n-n_3)^+$, which is
$$
\frac{t^{n_3+1} [ (n_3+1)^2 t^2 - 2(n_3^2-2) + (n_3-1)^2 ]}{(1-t)^4}.
$$

\myexam{A Small Problem of Queens and Bishops, continued from Example \ref{X:small}}\label{X:gfsmall}
The generating functions are, for bishops, 
\[
N_{B_2}(t) = \frac{A_2(t)}{(1-t)^3} - 2\frac{t^2}{(1-t)^2} = \frac{2t^3 - t^2 + t}{(1-t)^3} ,
\]
and for queens, 
\[
N_{Q_2}(t) = \frac{2t^3-t^2+t}{(1-t)^3} - \frac{t}{(1-t)^2} = \frac{2t^3}{(1-t)^3} .
\]
\end{exam}

%%%%%%%%%%%%%%
\section{Some chess pieces, real and fairy}

Let us look at some more examples.  We assume in all of them that $m \geq 1$ and there is one piece in each row (all $q_j=1$).  We write $\nu_m(n)$ for the number of nonattacking arrangements with $m$ rows, $N_m(t) := \sum_{n=0}^\infty \nu_m(n) t^n$ for its generating function, and $\bar\nu_m(n)$ for the polynomial when $n$ is large.  
In all the formulas, we assume $n\geq0$.  

Every example has $\nu_1(n) = n$ for $n\geq0$, and $N_1(t) = t/{(1-t)^{m+1}}.$  Also, all vertex weights are zero in the weighted integral gain graph for each piece.

What we would most like would be to find a pattern in the generating function numerator for a fixed piece as the height $m$ varies; but one would expect any regularity to show up only when the height is somewhat large because when it is small the moves are strongly constrained by the narrowness of the strip.  The solutions we calculate are, unfortunately, too small to show a pattern. 

We calculated the function $\chiz_\Phi(n)$ in several examples by reducing to gain graphs whose edges are all loops (and consequently may be discarded, unless a loop gain is zero, in which case the graph itself may be discarded), using deletion and contraction, as illustrated in Examples \ref{X:small} and \ref{X:gfsmall}.  The relationship between the graph quantity $\chiz_\Phi$ and $\nu_m$ is that they are equal, according to Equation \eqref{E:convert}.  We were able to do hand calculations for bishops and queens with $m\leq3$ and knights with $m\leq4$; we verified and extended these results by a computer package that carries out deletion and contraction of weighted integral gain graphs \cite{WIGGprogram} with the appropriate simplifications; that is, it discards nonzero loops and graphs with zero loops.

\myexam{Rooks}
 Rooks, of course, are classically easy: $\nu_m(n) = (n)_m$, the falling factorial.  The theorem says this is a polynomial for $n>0$, which is true.  
The generating function is 
$N_m(t) = \sum_{n=0}^\infty (n)_m t^n = {t^m m!}/{(1-t)^{m+1}},$ 
as is well known.
\end{exam}

\myexam{Bishops} \label{Bishops}
Proposition \ref{P:exact} applies with $b=1$.  Thus, $\nu_m(n)$ is a polynomial for $n \geq \frac12 (m^2 - 3)$ but not if we extend the domain by one more integer, to $n \geq \frac12 (m^2 - 5)$.

The bishops gain graph $B_m$ has edges $\mu v_iv_j$ with gains $\pm(j-i)$ for $i<j$.  We found that
\begin{align*} 
\nu_2(n) &= n^2 - 2(n-1)^+ , \\
\nu_3(n) &= n^3 - n\big\{ 4(n-1)^+ + 2(n-2)^+ \big\} \\
  &\quad + \big\{ 2(n-1)^+ + 4(n-2)^+ + 4(n-3)^+ \big\} , \\
\nu_4(n) &= n^4 - n^2\big\{ 6(n-1)^+ + 4(n-2)^+ + 2(n-3)^+ \big\} \\
  &\quad + n\big\{ 4(n-1)^+ + 12(n-2)^+ + 16(n-3)^+ + 4(n-4)^+ + 4(n-5)^+\big\} \\
  &\quad + \big\{ 4\big[(n-1)^+\big]^2 + 4(n-1)^+(n-3)^+ + 4\big[(n-2)^+\big]^2 \big\} \\
  &\quad - \big\{ 2(n-1)^+ + 8(n-2)^+ + 34(n-3)^+ + 12(n-4)^+ \\
  &\quad\qquad  + 20(n-5)^+ + 4(n-6)^+ + 2(n-7)^+\big\} .
\end{align*} 
The polynomials for large $n$ are
\begin{align*}
\bar\nu_2(n) &= n^2 - 2n + 2			&\text{ for } n \geq 1 , \\
\bar\nu_3(n) &= n^3 - 6n^2 + 18n - 22		&\text{ for } n \geq 3 , \\
\bar\nu_4(n) &= n^4 - 12n^3 + 72n^2 - 234n + 338	&\text{ for } n \geq 7  , \\
\bar\nu_5(n) &= n^5 - 20n^4 + 200n^3 - 1192n^2 + 4132n - 6562	&\text{ for } n \geq 11 , \\
\bar\nu_6(n) &= n^6 - 30n^5 + 450n^4 - 4198n^3 + 25238n^2 - 91572n + 155220	&\text{ for } n \geq 17 .
\end{align*}
The lower limits of validity can be read off from the piecewise-polynomial formulas (shown for $m \leq 4$; we computed $m=5,6$ as well) and are as predicted by Proposition \ref{P:exact}.  
The coefficients of $n^{m-1}$ and $n^{m-2}$ in the polynomials $\bar\nu_m(n)$ are the values $-2\binom m2$ and $2\binom{m}{2}^2$.  The first can be explained intuitively by the fact that when $n$ is large a bishop in row $i$ blocks $2(m-i)$ squares in higher rows (see Proposition \ref{P:secondterm}).  We have no explanation of the second.

Generating functions are
\begin{align*}%\label{E:}
N_2(t) &= \frac{2t^3-t^2+t}{(1-t)^3}, \qquad
N_3(t) = \frac{4t^6-4t^5-2t^4-t^3+22t^2-13t}{(1-t)^4}, \\
N_4(t) &= \frac{2t^{11}-2t^{10}+14t^9-38t^8+54t^7-74t^6+64t^5-43t^4+213t^3-247t^2+81t}{(1-t)^5}.
\end{align*}

\begin{figure}[thb]
\includegraphics[scale=.8]{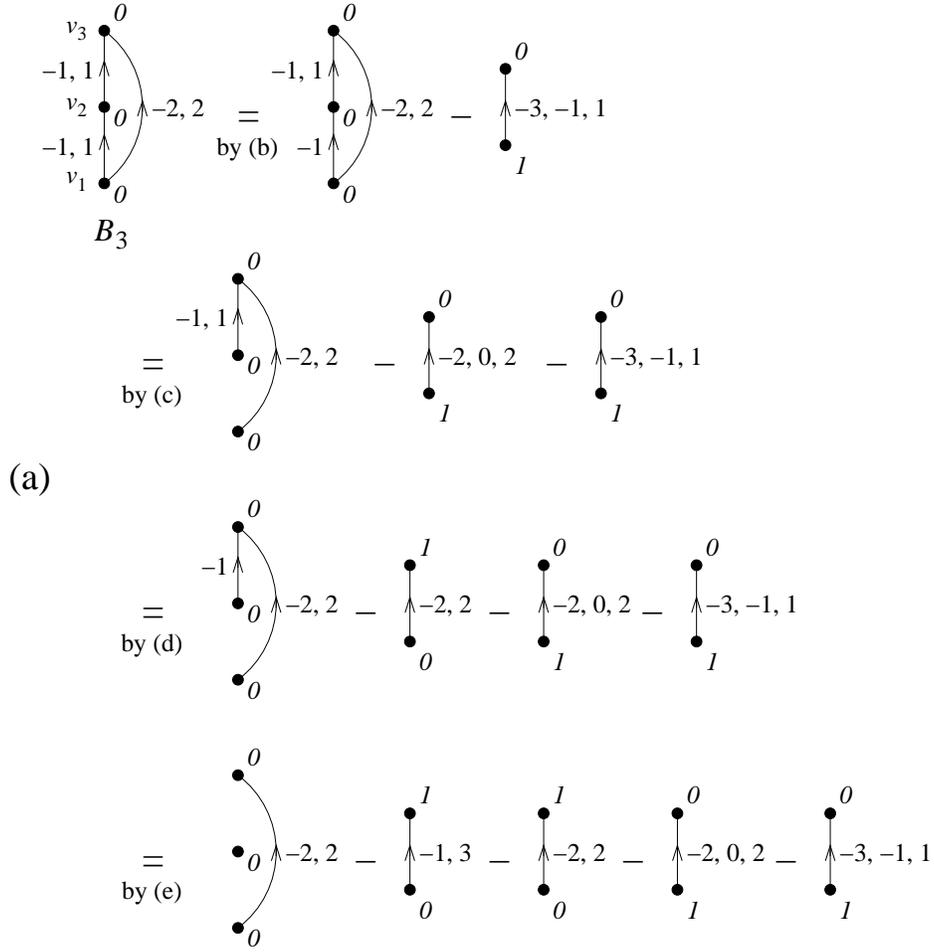}
\caption{(a)  The weighted integral gain graph for $3$ bishops in $3$ rows, showing the contraction process for (b) the first, (c) the second, (d) the third, and (e) fourth steps of deletion and contraction.  $e$ denotes $1v_1v_2$, $e' = -1v_1v_2$, and $e'' = 1v_2v_3$.  For simplicity, parallel edges are drawn as a single edge with multiple gains.  Nonzero loops are omitted.}
\label{F:b3}
\end{figure}

As a further illustration of the method of gain graphs we give some details of the solution for $\nu_3$.  The gain graph $B_3$, shown in Figure \ref{F:b3}(a), has vertices $v_1, v_2, v_3$ with $v_iv_j$-edges having gains $\pm(j-i)$ for $i<j$.  
We delete and contract edges in the arbitrarily chosen order $1v_1v_2$, $-1v_1v_2$, $1v_2v_3$, $-1v_2v_3$, after which each component is a vertex or a multiple edge.  The treatment of these multiple edges is as in Example \ref{X:multiple}, so is not shown here.  

\begin{figure}[thb]
\includegraphics[scale=.8]{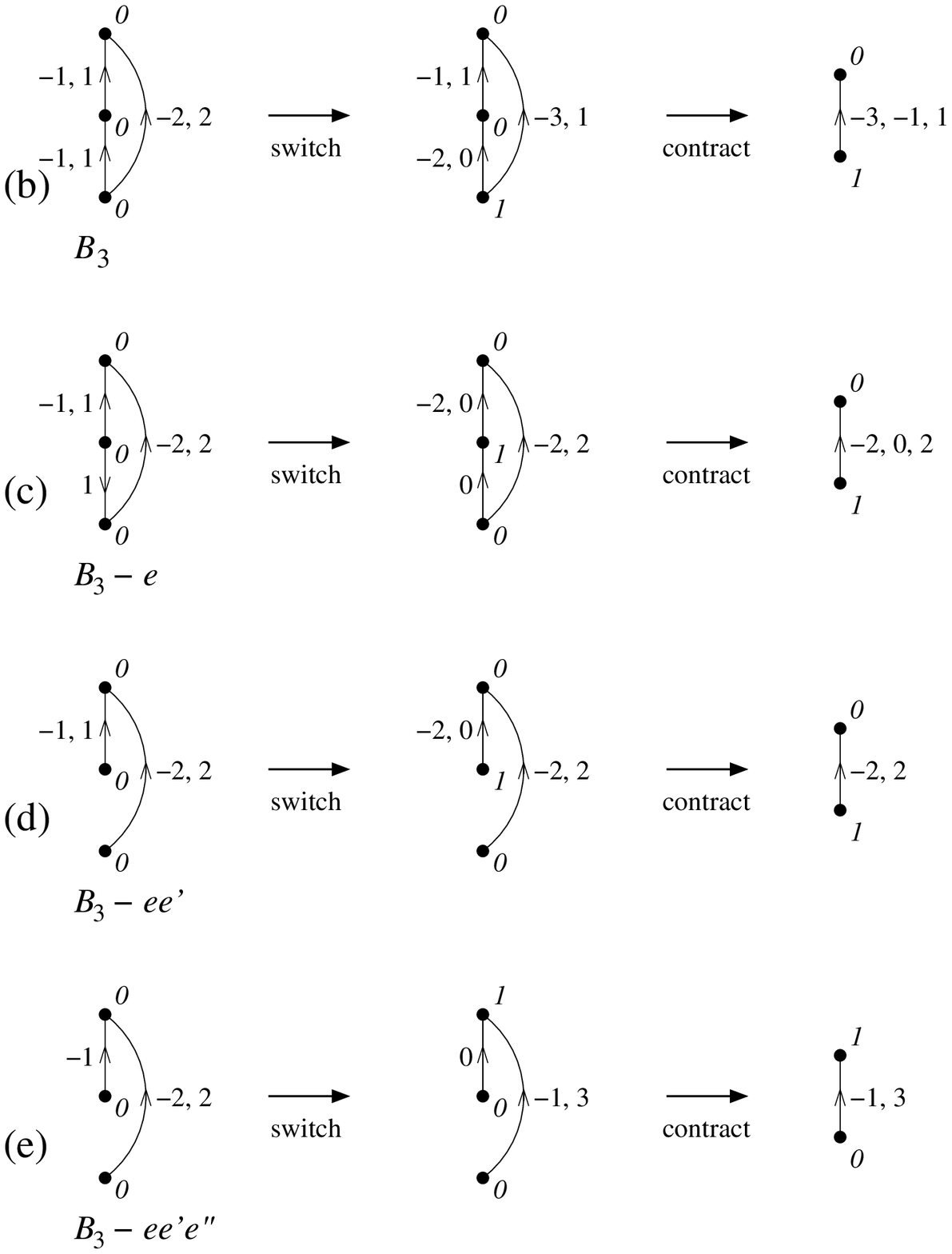}
%
%\caption{Parts (b--e) of Figure 3.}
\label{F:b3be}
\end{figure}
\end{exam}

\myexam{Queens} \label{Queens}
As with bishops, Proposition \ref{P:exact} applies with $b=1$.  Thus, $\nu_m(n)$ is a polynomial for $n \geq \frac12 (m^2 - 3)$ but not for $n \geq \frac12 (m^2 - 5)$.

The queens gain graph $Q_m$ has edges $\mu v_iv_j$ with gains $\mu = 0, \pm(j-i)$ for $i<j$; it is the bishops graph $B_m$ together with $0$ edges.  The results are:
\begin{align*} 
\nu_2(n) &= n^2 - \big\{ n + 2(n-1)^+ \big\} , \\
\nu_3(n) &= n^3 - n \big\{ 3n + 4(n-1)^+ + 2(n-2)^+ \big\} \\
  &\quad + \big\{ 2n + 8(n-1)^+ + 8(n-2)^+ + 4(n-3)^+ \big\} , \\
\nu_4(n) &= n^4 - n^2\big\{ 6n + 6(n-1)^+ + 4(n-2)^+ + 2(n-3)^+ \big\} \\
  &\quad + n\big\{ 11n + 30(n-1)^+ + 32(n-2)^+ + 26(n-3)^+ + 4(n-4)^+ + 4(n-5)^+\big\} \\
  &\quad + \big\{ 4\big[(n-1)^+\big]^2 + 4(n-1)^+(n-3)^+ + 4\big[(n-2)^+\big]^2 \big\} \\
  &\quad - \big\{ 6n + 40(n-1)^+ + 72(n-2)^+ + 94(n-3)^+ + 36(n-4)^+ \\
  &\quad\qquad + 32(n-5)^+ + 4(n-6)^+ + 2(n-7)^+\big\} .
\end{align*} 
The polynomials for large $n$ are
\begin{align*}
\bar\nu_2(n) &= n^2 - 3n + 2			&\text{ for } n \geq 1 , \\
\bar\nu_3(n) &= n^3 - 9n^2 + 30n - 36		&\text{ for } n \geq 3 , \\
\bar\nu_4(n) &= n^4 - 18n^3 + 139n^2 - 534n + 840	&\text{ for } n \geq 7 , \\
\bar\nu_5(n) &= n^5 - 30n^4 + 407n^3 - 3098n^2 + 13104n - 24332	&\text{ for } n \geq 11 , \\
\bar\nu_6(n) &= n^6 - 45n^5 + 943n^4 - 11755n^3 + 91480n^2 - 418390n + 870920	&\text{ for } n \geq 17 .
\end{align*}
The ranges of validity of the polynomials agree with Proposition \ref{P:exact}.  The second term in the polynomial has coefficient $-3\binom{m}{2}$, which as with bishops is the sum of the number of squares excluded by each piece in rows higher than itself.  The third coefficient does not follow an easily identifiable pattern.

The generating functions of the exact numbers are
\begin{align*}%\label{E:}
N_2(t) &= \frac{2t^3}{(1-t)^3} , \qquad
N_3(t) = \frac{4t^6-4t^4-4t^3+24t^2-14t}{(1-t)^4}, \\
N_4(t) &= \frac{2t^{11}-2t^{10}+26t^9-50t^8+78t^7-130t^6+66t^5-42t^4+484t^3-640t^2+232t}{(1-t)^5} .
\end{align*}
\end{exam}

\myexam{Knights}\label{Knights}
A widest knight's path lets the knight rise from the bottom row to the top, one row at a time, moving 2 steps to the right at each row.  There are $m-1$ moves that cover a width of $2m-2$.  We conclude that $\nu_m(n)$ is a polynomial for $n\geq 2m-2$ but not for $n\geq 2m-3$. 

The knights gain graph $Kt_m$ has edges $\pm2 v_iv_{i+1}$ and $\pm1 v_iv_{i+2}$.

We calculated the first few cases, $m\leq3$ by hand and $m\leq6$ on the computer, via gain graphs.
\begin{align*} 
\nu_2(n) &= n^2 - 2(n-2)^+ , \\
\nu_3(n) &= n^3 - n \big\{ 2(n-1)^+ + 4(n-2)^+ \big\} + \big\{ 6(n-2)^+ + 4(n-3)^+ + 2(n-4)^+ \big\} , \\
\nu_4(n) &= n^4 - \big\{ 4n^2(n-1)^+ + 6n^2(n-2)^+ \big\} \\
&\quad + \big\{ 16n(n-2)^+ + 12n(n-3)^+ + 4n(n-4)^+ + 4[(n-1)^+]^2 + 4[(n-2)^+]^2 \big\} \\ 
&\quad - \big\{ 12(n-2)^+ + 24(n-3)^+ + 18(n-4)^+ + 4(n-5)^+ + 2(n-6)^+ \big\} .
\end{align*} 
The polynomials for large $n$ are
\begin{align*}
\bar\nu_2(n) &= n^2 - 2n + 4			&\text{ for } n \geq 2 , \\
\bar\nu_3(n) &= n^3 - 6n^2 + 22n - 32		&\text{ for } n \geq 4 , \\
\bar\nu_4(n) &= n^4 - 10n^3 + 56n^2 - 168n + 220	&\text{ for } n \geq 6 , \\
\bar\nu_5(n) &= n^5 - 14n^4 + 106n^3 - 478n^2 + 1248n - 1480	&\text{ for } n \geq 8 , \\
\bar\nu_6(n) &= n^6 - 18n^5 + 172n^4 - 1028n^3 + 3956n^2 - 9154n + 9852	&\text{ for } n \geq 10.
\end{align*}
The ranges of validity, read off from the piecewise-polynomial formulas, agree with the preceding calculation.  
The size of the second coefficient, as with bishops and queens, is $c_1 =  4m-6$, which is the number of squares excluded in higher rows as well-separated knights are inserted from the bottom row up, that is, 4 squares for each knight except the one in the next-to-last row, which excludes 2 squares in the last row, and the last knight, which excludes no squares since there are no rows below it.  
(See Proposition \ref{P:secondterm}.)  
The third coefficient appears to be a quadratic function of $m$, $c_2=8m^2-22m+16$. (As with bishops, we have not proved this.)

As for generating functions, we have 
\begin{align*}%\label{E:}
N_2(t) &= \frac{2t^4-2t^3+t^2+t}{(1-t)^3} , \qquad
N_3(t) = \frac{2t^7-8t^4+3t^3+24t^2-15t}{(1-t)^4} , \\
N_4(t) &= \frac{2t^{10}-2t^9+12t^8-20t^7-6t^6+18t^5-25t^4+181t^3-199t^2+63t}{(1-t)^5} .
\end{align*}
\end{exam}

\myexam{Nightriders} \label{Nightriders}
This fairy chess piece has the move of a knight extended indefinitely in a straight line; that is, it moves along lines of slope $\pm 2^{\pm1}$.  Therefore Proposition \ref{P:exact} applies with $b=2$, and $\nu_m(n)$ is a polynomial for $n\geq 2 \lfloor m^2/2 \rfloor - 2$ but not for $n\geq 2 \lfloor m^2/2 \rfloor - 3$.

The gain graph $Nr_m$ for nightriders has edges $\pm2(j-i) v_iv_j$ if $j-i$ is nonzero, and edges $\pm\frac12(j-i) v_iv_j$ if $j-i$ is nonzero and even.  The former correspond to moves of any number of rows in the directions $\pm(\pm2,1)$ and the latter to moves of an even number of rows in the steeper directions $\pm(\pm1,2)$.

Using the computer, we calculated the same data as in the previous examples for $m\leq 6$.  The counting functions for $m=1,2$ are the same as for the knight.  For $2\leq m \leq 4,$ they are
\begin{align*} 
\nu_2(n) &= n^2 - 2(n-2)^+ , \\
\nu_3(n) &= n^3 - n \big\{ 2(n-1)^+ + 6(n-2)^+ \big\} + \big\{ 10(n-2)^+ + 4(n-3)^+ + 2(n-4)^+ \big\} , \\
\nu_4(n) &= n^4 - \big\{ 4n^2(n-1)^+ + 12n^2(n-2)^+ \big\} \\
&\quad + \big\{ 40n(n-2)^+ + 16n(n-3)^+ + 24n(n-4)^+ \\ & \qquad\qquad + 4[(n-1)^+]^2 + 8(n-1)^+(n-2)^+ + 12[(n-2)^+]^2 \big\} \\ 
&\quad - \big\{ 58(n-2)^+ + 56(n-3)^+ + 100(n-4)^+ + 24(n-5)^+ + 24(n-6)^+ \big\} .
\end{align*} 
The polynomials for large $n$ are
\begin{align*}
\bar\nu_2(n) &= n^2 - 2n + 4			&\text{ for } n \geq 2 , \\
\bar\nu_3(n) &= n^3 - 8n^2 + 34n - 56		&\text{ for } n \geq 4 , \\
\bar\nu_4(n) &= n^4 - 16n^3 + 132n^2 - 566n + 1016	&\text{ for } n \geq 6 , \\
\bar\nu_5(n) &= n^5 - 28n^4 + 390n^3 - 3100n^2 + 13600n - 25676	&\text{ for } n \geq 8 , \\
\bar\nu_6(n) &= n^6 - 42n^5 + 876n^4 - 10974n^3 + 84720n^2 - 374678n + 730408	&\text{ for } n \geq 10 .
\end{align*}
The second coefficient, up to sign, is $2\binom{m}{2}+2\binom{\lfloor(m+1)/2\rfloor}{2}+2\binom{\lfloor m/2\rfloor}{2}$, as Proposition \ref{P:secondterm} implies.  
For nightriders, as with queens, the third coefficient $c_2$ does not have a discernable pattern.

The generating functions are 
\begin{align*}%\label{E:}
N_2(t) &= \frac{2t^4-2t^3+t^2+t}{(1-t)^3} , \qquad
N_3(t) = \frac{6t^7-8t^6+8t^5-16t^4+5t^3+32t^2-21t}{(1-t)^4} , \\
N_4(t) &= \frac{24t^{10}-48t^9+100t^8-196t^7+178t^6-114t^5+15t^4+144t^3-585t^2+209t}{(1-t)^5} .
\end{align*}
\end{exam}

%%%%%%%%%%%%%%
\section{The polynomial}\label{polynomial}

We observed in the examples that the polynomial $\bar\nu_m(n)$ which gives $\nu_m(n)$ for wide boards, i.e., for large values of $n$, has second term $-c_1 n^{m-1}$ where $c_1$ is the sum over all rows $i$ of the number of squares attacked in higher rows by a piece in row $i$.  There is a general result here.  Suppose there are $q_i$ pieces in row $i$.  The polynomial can be written in the form 
$$
\bar\nu_m(n) = \frac{n^q - c_1 n^{q-1} + c_2 n^{q-2} + \cdots}{q_1! \cdots q_m!}
$$ 
with $q = q_1 + \cdots + q_m$.

\begin{prop}\label{P:secondterm}
Assume that $(0,0)$ is a move.  Then the magnitude of the second coefficient of $\bar\nu_m(n)$ is 
$$
c_1 = \sum_{i=1}^m \left[ \binom{q_i}{2}+ q_i \sum_{j=i+1}^m q_j a_{ij} \right],
$$
where $a_{ij}$ is the number of squares in row $j$ that are attacked by a piece in row $i$.  
When every row has one piece, 
$$
c_1 = \sum_{i=1}^m  a_i ,
$$
where $a_i$ is the number of squares in rows $i+1,\ldots,m$ that are attacked by a piece in row $i$; thus, $c_1$ equals half the maximum total number of squares attacked (nonhorizontally) by one piece placed in each of the $m$ rows.
\end{prop}

One calculates the number of squares attacked ignoring any limitations due to inadequate board width.

\begin{proof}
The proof is by inclusion and exclusion using labelled pieces.  
Let $w$ be the greatest width of a nonhorizontal move that is possible on a board of height $m$, and assume $n$ is very large compared to $w$.  
Number the pieces $P_1,P_2,\ldots,P_q$ so that the pieces in higher rows have higher numbers.  

Imagine the board extended indefinitely to left and right by ``imaginary'' squares supplementing the ``real'' squares of the $m \times n$ strip.  Let $U$ be the set of all placements of the $q$ pieces, attacking or not, and let $A_{kl}$ be the set of placements in which $P_k$ attacks $P_l$ and the latter is in a higher row.  Let $i$ be the row of $P_k$ and $j$ that of $P_l$.  
Then 
$$
|U| = \prod_{h=1}^m (n)_{q_h} 
 = n^q - \sum_{h=1}^m \binom{q_h}{2} n^{q-1} + \text{ terms of lower degree} ,
$$
where $(n)_r$ is the falling factorial.

To calculate $|A_{kl}|$ we place all pieces except $P_l$, and then we put $P_l$ into any square (in row $j$) attacked by $P_k$, including occupied and imaginary squares.  Call the set of all such arrangements, in which $P_k$ attacks $P_l$ and the other pieces may be attacking, $A'$.  Some of these arrangements have $P_l$ in an imaginary square; call the set of those arrangements $B$.  Some of the arrangements have $P_l$ in an occupied square; call the set of those arrangements $C$.  Since $B$ and $C$ are disjoint, $|A_{kl}| = |A'| - |B| - |C|$.  We estimate $|A'|$, $|B|$, and $|C|$ in order to estimate $|A|$.

First, $|A'|$ is the number of ways to place $q_h$ pieces in each row $h$, except that only $q_i-1$ pieces appear in row $i$, and then to put $P_l$ in any of the $a_{ij}$ real or imaginary squares in row $j$ that are attacked by $P_k$.  Thus, 
\begin{align*}%\label{E:}
|A'| &=  (n)_{q_1} \cdots (n)_{q_{j-1}} (n)_{q_j-1} (n)_{q_{j+1}} \cdots (n)_{q_m}  a_{ij} 
 = a_{ij} n^{q-1} + \text{ lower terms} .
\end{align*}

Second, $|B|$ counts only placements where $P_k$ lies within distance $w$ from the ends of the strip, because it must be able to attack an imaginary square.  If we count such arrangements by placing $P_k$ and $P_l$ last, there are no more than $n^{q-2}$ ways to place the other pieces in real squares and then at most $2w$ ways to lay down $P_k$ and $a_{ij}$ ways to put down $P_l$.  That is, $|B|$ is bounded by a polynomial in $n$ of degree less than $q-1$.

Third, $|C|$ counts placements where $P_l$ is attacked by $P_k$ and overlies another piece in its row.  Let us again lay down $P_k$ and $P_l$ last.  The number of places for $P_k$ is not more than $a_{ij}(q_j-1)$, because $P_k$ has to be attacked by one of the $q_j-1$ pieces already placed in row $j$,  and because $a_{ji} = a_{ij}$ by symmetry of moves.  Then the number of places for $P_l$ is also no more than $q_j-1$.  Thus, $|C|$ is bounded by a polynomial in $n$ of degree less than $q-1$.

It follows that $|A_{kl}| = a_{ij} n^{q-1}$ + a quantity bounded by a polynomial of degree at most $n^{q-2}$.

The size of the intersection of two sets of bad placements, $|A_{kl} \cap A_{k'l'}|$, is bounded by the number of placements in which $P_l$ and $P_{l'}$ are attacked by $P_k$ and $P_{k'}$, respectively, which is in turn bounded by $a_{ij}a_{i'j'}n^{q-2}$, a polynomial of degree less than $q-1$.  
Likewise, the intersection of more than two sets $A_{kl}$ is bounded in size by a polynomial of degree less than $q-1$.  

The proposition now follows from the formula of inclusion and exclusion combined with the fact that $\bar\nu_m(n)$ is known to be a polynomial of degree $q$.  Taking only the terms of highest order, and letting $(k,l)$ range over pairs of piece indices for which $P_k$ is in row $i$ and $P_l$ is in row $j$ with $i<j$, 
\begin{align*}%\label{E:}
\Big| \bigcap_{k,l} \bar A_{kl} \Big| &= |U| - \sum_{k,l} |A_{kl}| + \text{ lower terms} \\
 &= \Big[ n^q - \sum_{h=1}^m \binom{q_h}{2} n^{q-1} \Big] 
 - \sum_{i<j\leq m}  q_i q_j a_{ij} \, n^{q-1}  + \text{ lower terms} .
\qedhere
\end{align*}
\end{proof}

The especial importance of $c_1$ is shown by a probabilistic application.\footnote{We thank undergraduate student Christian Noack for raising this question.}

\begin{cor}\label{C:probability}
The probability that a random assignment of one piece in each row is nonattacking is asymptotic to $1 - c_1 n\inv$.  Assuming that $(0,0)$ is a move, the probability that a random assignment of $q_i$ pieces in distinct positions in row $i$ is nonattacking is asymptotic to $1 - n\inv {\sum\sum}_{i<j} \, q_i q_j a_{ij}.$
\end{cor}

\begin{proof}
The first probability is $\nu_m(n)/n^m = 1 - c_1 n\inv + O(n^{-2})$.

The second probability is 
\begin{align*}
\frac{\nu_m(n)}{\prod_{h=1}^m \binom{n}{q_h}} 
&= \frac{ 1 - c_1 n\inv + O(n^{-2}) }{ 1 - n\inv \sum_{i=1}^m \binom{q_i}{2} + O(n^{-2}) } \\
&= \big[1 - c_1 n\inv + O(n^{-2}) \big] \Big[ 1 + n\inv \sum_{i=1}^m \binom{q_i}{2} + O(n^{-2}) \Big]
\intertext{by the geometric series (when $n$ is sufficiently large),}
&= 1 - n\inv \Big[ c_1 - \sum_{i=1}^m \binom{q_i}{2} \Big] + O(n^{-2}) .  
\qedhere
\end{align*}
\end{proof}

We did not find a formula for $c_2$.  The solution for bishops and the proof of Proposition \ref{P:secondterm} suggest there may be one but that it is not usually simple.

%%%%%%%%%%%%%%%%%%%%
\section*{Acknowledgements}

We thank the referees for their well-judged suggestions towards improving the exposition.

\emph{Note added February, 2010.}  
After acceptance of this article, V\'aclav Kot{\v e}{\v s}ovec (January 25, 2010) informed us of the prior discovery of the polynomials $\bar\nu_k(n)$ for queens (Example \ref{Queens}) by Pauls \cite{Pauls} ($k=3$, 1874), Tarry \cite{Tarry} ($k=4$, 1890) and Kot{\v e}{\v s}ovec \cite{Azemard, K1996} ($k=5,6,7$, 1992).  Initial values of the polynomials can be found in the On-Line Encyclopedia of Integer Sequences \cite{OEIS}, sequences A061989--A061993.  
The references were provided by Kot{\v e}{\v s}ovec.
Kot{\v e}{\v s}ovec's Web site \cite{Kotesovec} has many other formulas for placing non-attacking chess pieces.  Furthermore, Kot{\v e}{\v s}ovec has produced (February, 2010) a conjectural simple formula for the third coefficient of the queens polynomials.

%%%%%%%%%%%%%%%%%%%%

\end{document}